\documentclass[leqno,12pt]{article}
\usepackage{amsmath, amssymb}
\usepackage{amsthm}

\theoremstyle{plain}
\newtheorem{theorem}{\indent\sc Theorem}[section]
\newtheorem{lemma}[theorem]{\indent\sc Lemma}
\newtheorem{corollary}[theorem]{\indent\sc Corollary}

\theoremstyle{definition}
\newtheorem{definition}[theorem]{\indent\sc Definition}
\newtheorem{remark}[theorem]{\indent\sc Remark}



\newcommand{\field}[1]{\mathbb{#1}}

\newcommand{\R}{\field{R}}

\newcommand{\beeq}{\begin{equation}}
\newcommand{\eneq}{\end{equation}}
\newcommand{\la}{\label}
\newcommand{\pref}[1]{(\ref{#1})}
\newcommand{\ba}{\begin{array}}
\newcommand{\ea}{\end{array}}

\newcommand{\ds}{\displaystyle}
\newcommand{\be}{\begin{equation}}
\newcommand{\ee}{\end{equation}}
\newcommand{\bea}{\begin{eqnarray}}
\newcommand{\eea}{\end{eqnarray}}
\newcommand{\nn}{\nonumber}
\newcommand{\zume}{\!\!\!}

\newcommand{\intmean}{{\int\hspace{-12.5pt}-}}

\newcommand{\A}{{\cal A}}

\newcommand{\supp}{\mathrm{supp}}
\begin{document}
{\centerline{\bf{ Regularity of minimizers of some variational integrals\\ with discontinuity}}
%

\center{
Maria Alessandra Ragusa\\
Dipartimento di Matematica e Informatica, Universit\`{a} di Catania,
Viale Andrea Doria, 6-95125 Catania, Italy,\\
e-mail:maragusa@dmi.unict.it
}
\center{
Atsushi Tachikawa\\
Department of Mathematics, Faculty of Science and Technology,
Tokyo University of Science, Noda, Chiba, 278-8510, Japan\\
e-mail:tachikawa$\_$atsushi@ma.noda.tus.ac.jp 
}

\abstract{
We prove regularity properties in the vector valued case for minimizers
of variational integrals of the form $$  \A(u) = \int_\Omega A(x,u,Du) dx$$
where the integrand $A(x,u,Du)$ is not necessarily continuous
respect to the
variable $x,$ grows polinomially like $|\xi|^p,$
$p \geq 2.$
}
$$\,$$
Keywoards: variational problem, minimizer, partial regularity 
MSC: 35J10, 35B65, 35N10, 
46E30, 35R05
\section{Introduction}

In this note we consider the regularity problem of
minimizers of the variational integral
\begin{equation}\label{def-A}
\A(u) = \int_\Omega A(x,u,Du) dx
\end{equation} where $\Omega$ is a bounded
domain of $\R^m,$ $u:\Omega \to \R^n$ is a
mapping in a suitable Sobolev space,
$Du = (D_\alpha u^i)~~(\alpha =1,\ldots,m,~i=1,\ldots,n)$.
The nonnegative integrand function
$A: \Omega \times \R^n \times \R^{mn}$ is in the class \it VMO \rm
with respect to the variable $x,$ continuous in $u$
and of class $C^2$ with respect to $Du.$
It is also assumed that for some $p\geq 2$
there exist two
constants $\lambda_1$  and $\Lambda_1$ such that
\begin{equation}\label{def-B}
     \lambda_1 (
     1+|\xi|^p) \leq A(x,u,\xi) \leq \Lambda_1(
     1 + |\xi|^p), \quad \forall (x,u,\xi)\in \Omega\times \R^n \times
     \R^{mn},
\end{equation}
A minimizer for the functional $\A$ is a function $u \in
W^{1,p}(\Omega , \R^n)$ such that for every
$\varphi \in W^{1,p}_0(\Omega , \R^n)$
$$
\A(u;\supp \varphi) \leq \A(u+\varphi; \supp \varphi).
$$

For the case that $A(x,u,\xi)$ is continuous in $x$,
many sharp regularity results for minimizers of $\cal{A}$
have been already known. (eg. \cite{Gia-book1, Gia-book2, GL, Giusti-book}.)
On the other hand, when $A(\cdot ,u,\xi)$ is assumed only to be
$L^\infty$,
we can not expect the regularity of minimizers in general,
as a famous example due to De Giorgi contained in \cite{DeG68} asserts.
So, it seems to be natural to consider the regularity problems
for $A(x,u,\xi)$ with \it ``mild'' \rm discontinuity with respect to $x$.
In 1996 Huang in \cite{QH} investigates regularity results for the elliptic
system
$$
- D_\alpha(a_{ij}^{\alpha \beta}(x) D_\beta u^j ) \,=\, g_i(x) - \mathrm{div}
f^i (x), \qquad i,j=1,\ldots,n; \alpha, \beta =1,\ldots,m
$$
assuming that $a_{ij}^{\alpha \beta}$
belong to the Sarason class $VMO$ of vanishing mean oscillation functions,
then he generalizes Acquistapace's \cite{A} and
Campanato's results \cite[p.88, Theorem 3.2]{Gia-book1}.
Campanato showed regularity properties under the assumption that the
coefficients $a_{ij}^{\alpha\beta}$
are in $C^\alpha (\Omega).$
Acquistapace refined the results by Campanato,
considering  
the coefficients in the class so-called  ``small multipliers of $BMO$''.

In the recent study made by Dan\v{e}\v{c}ek and Viszus \cite{DV}
they consider the following functional
$$
\int_\Omega \left\{  A^{\alpha\beta}_{ij}(x) D_\alpha u^i D_\beta u^j + g(x,u,Du)\right\} dx,
$$
where $A^{\alpha\beta}_{ij}$ are in general discontinuous, more
precisely belong to the vanishing mean oscillation class ($VMO$ class)
and satisfy strong ellipticity condition while the lower order term $g$
is a Charath\'eodory function satisfying the following growth condition
\begin{equation}\label{(6)}
|g(x,u,z)|\,\leq \, f(x) + H |z|^\kappa
\end{equation}
where $f \geq 0, $ a.e. in $\Omega,$ $f \in L^p(\Omega),$ $2<p\leq \infty
,$ $H\geq 0,$ $0\leq \kappa < 2.$

We also recall the paper by Di Gironimo, Esposito and Sgambati \cite{DES}
where is treated the Morrey regularity for minimizers of the functional
\[
\int_\Omega A^{\alpha\beta}_{ij}(x,u) D_\alpha u^i D_\beta u^j dx,
\]
where $(A^{\alpha\beta}_{ij}(x,u))$ are elliptic and of the $VMO$ class
in the variable $x$.

In \cite{RT1} the authors extend the results of \cite{DV} and
\cite{DES} to the case that the functional is given by
$$
\int_\Omega \left\{  A^{\alpha\beta}_{ij}(x,u)
D_\alpha u^i D_\beta u^j + g(x,u,Du)\right\} dx.
$$
In the note \cite{RT2}, is studied the Morrey regularity for
minimizer of the more general functionals
\begin{equation}\label{(A)}
  \A(u) = \int_\Omega A(x,u,Du) dx,
\end{equation}
where $A(x,u,\xi)$ is a nonnegative function defined on
$\Omega \times \R^n \times \R^{mn}$ which is of class
$VMO$ as a function of $x$, continuous in $u$ and of class $C^2$
with respect to $\xi$.
We point out that is assumed  that for some positive constants
$\mu_0 \leq \mu_1$,
\begin{equation}\label{(2-growth)}
     \mu_0|\xi|^2 \leq A(x,u,\xi) \leq \mu_1 |\xi|^2
\quad \forall(x,u,\xi)\in \Omega\times \R^n \times \R^{mn}.
\end{equation}

We point out that
in the above mentioned papers concerning
functionals given by integrals with $VMO$ class
integrands,
we have considered quadratic growth functionals.

The super quadratic cases with continuous coefficients are treated in
\cite{CamAtti} and \cite{GM-Manus}.

In the present note we investigate the partial regularity of
the minima of $\A,$ defined by \eqref{def-A} under $p$-growth
hypothesis of the integrand function $A,$
$p\geq 2$.
This study can be considered as an improving of \cite{RT1} and
\cite{RT2} because of the growth condition is more general.

\section{Definitions and Preliminary Tools}\label{pre-res}
\setcounter{equation}{0}
\setcounter{theorem}{0}

In the sequel we set
\[
Q(x,R) = \{ y \in \R^m~\colon~ |y^\alpha - x^\alpha| <R,~ \alpha=
1,...,m\}
\]
a generic cube in $\R^m$ having center $x$ and side $2R$.

Let us now give some useful definitions, starting to the Morrey space
$L^{p,\lambda}$ .
\begin{definition}\label{def2.1}(see \cite{KJF}).
Let $1\leq p<\infty, 0 \leq \lambda <m.$ A measurable function
$ G \in L^p (\Omega, \R^n)$ belongs to the Morrey
class $L^{p,\lambda} (\Omega ,\R^n)$ if
\[
\|G \|_{L^{p,\lambda} (\Omega)} =
\sup_{\underset{x \in \Omega}{0<\rho < \mathrm{diam}\, \Omega }}
\frac{1}{\rho^\lambda}
\int\limits_{\Omega \cap Q(x,\rho)}
|G (y)|^p dy < + \infty
\]
where $Q(x,\rho)$ ranges in the class of  the
cubes of ${\R}^{m}$.
\end{definition}
\begin{definition}\label{def2.2}
Let $H \in L^1(\Omega , \R^n )$ we set the integral average
$H_{x,R}$ by
\[ H_{x,R}\,=\, \intmean_{\Omega \cap Q(x,R)} H (y)\, dy \,=\,
\frac{1}{|\Omega \cap Q(x,R)|} \int\limits_{\Omega \cap Q(x,R) }H (y) \,dy
\]
where $|\Omega \cap Q(x,R)|$ is the Lebesgue measure of $\Omega \cap
Q(x,R).$

In the case that we are not interested in specifying which the
center is, we simply write
$H_R. $
\end{definition}

Let us introduce the Bounded Mean Oscillation class.

\begin{definition}\label{def2.3}(see \cite{JN}).
Let
$H \in L^1_{\mathrm{loc}}({\R}^m).$ We say that
$H$ belongs
to $BMO({\R}^m)$ if
\[
\| H \|_* \equiv
\sup_{Q(x,R) } \frac{1}{|Q(x,R)|}
\int\limits_{Q(x,R)} |H(y)-H_{x,R}| dy < \infty .
\]
\end{definition}

Let us now introduce the space of vanishing mean oscillation functions
(\cite{Sarason}).

\begin{definition}\label{def2.4}
If $H \in BMO({\R}^m)$  and
\[
\eta (H;R) =
\sup_{\rho \leq R } \frac{1}{|Q(x,\rho)|}
\int\limits_{Q(x,\rho)} |H(y)-H_{\rho}| dy
\]
We have that $H\in VMO(\Omega)$ if
\[
\lim_{R \to 0} \eta (H; R) =0.
\]
\end{definition}

Throughout the present paper we consider $p \geq 2$ and $u:\Omega \to \R^n$
a minimizer of the functional
$$
\A(u) = \int_\Omega A(x,u,Du) dx
$$
where the hypothesis on the integrand function $A(x,u,\xi)$ are the
following.

%
\begin{enumerate}
\item[(A-1)]  For every $(u,\xi)\in \R^n \times \R^{mn}$,
    $A(\cdot, u, \xi) \in VMO(\Omega)$ and the mean oscillation of
    $A(\cdot, u, \xi)/|\xi|^p $ vanishes uniformly with respect to
    $u, \xi$ in the following sense:
    there exist a positive number $\rho_0$ and a function
    $\sigma(z, \rho) : \R^m \times [0, \rho_0) \to [0,\infty)$ with
    \begin{equation}
    \lim_{R\to 0} \sup_{\rho <R}
    \intmean_{Q(0,\rho)\cap \Omega} \sigma(z,\rho) dz =0,
    \end{equation}
    such that
    $A(\cdot , u ,\xi)$ satisfies
    for every $x\in \overline{\Omega}$ and
    $y \in Q(x,\rho_0)\cap \Omega$
    \begin{equation}\la{A-1}
	    \big|A(y,u,\xi) -
        A_{x,\rho}(u,\xi)
        \big| \leq \sigma(x-y,\rho)(
        1 + |\xi|^2)^{\frac{ p}{2}}
          ~~~~ \forall (u,\xi)\in \R^n \times \R^{mn},
    \end{equation}
    where
    \[
        A_{x,\rho}(u,\xi)=\intmean_{Q(x,\rho)\cap \Omega} A(y,u,\xi) dy.
    \]
\item[(A-2)]  For every $x\in \Omega$, $\xi \in \R^{mn}$ and $u,v\in \R^n$,
    \[ \big|A(x,u,\xi)-A(x,v,\xi)\big|
       \leq (
       1 + |\xi|^2 )^{\frac{p}{2}}
       \omega(|u-v|^2)
    \]
  where $\omega$ is some monotone increasing concave function
  with $\omega(0)  =0$,
\item[(A-3)] For almost all $x\in \Omega$ and all $u\in \R^n$,
    $A(x,u,\cdot) \in C^2(\R^{mn})$,
\item[(A-4)] There exist positive constants $\lambda_1,$ $\Lambda_1$ such that
     \[
     \lambda_1 (
     1+|\xi|^p) \leq A(x,u,\xi) \leq \Lambda_1(
     1 + |\xi|^p)
     \]
     \[
           \lambda_1 (1+|\eta|^p) \leq
           \frac{\partial^2 A(x,u,\xi)}{\partial \xi^i_\alpha
         \partial \xi^j_\beta} \eta^i_\alpha \eta^j_\beta \leq \Lambda_1(
           1 + |\eta|^p)
         \]
     for all $(x,u,\xi,\eta)\in \Omega\times \R^n \times \R^{mn}
\times \R^{mn}$.
\end{enumerate}

Let us state the main theorem of the paper concerning the partial
regularity of the minimizers of the functionals $\A.$

\begin{theorem}\la{mainTh}
Assume that $\Omega\subset \R^m$ is a bounded domain with
sufficiently smooth boundary $\partial \Omega$ and that $p\geq 2$.
Let $u \in
H^{1,p}(\Omega, \R^n)\,$  a
minimizer of the functional
\[
\A(u,\Omega) = \int_\Omega A(x,u,Du) dx
\]
in the class
\[
X_g (\Omega) = \{u \in H^{1,p}(\Omega)~;~
u-g \in H^{1,p}_0(\Omega)\}
\]
for a given boundary data $g \in H^{1,s}(\Omega)$
with
$s>p$.
Suppose that assumptions \rm (A-1), (A-2), (A-3)\it and \rm (A-4) \it
are satisfied.
Then, for some positive $\varepsilon$, for every
$0<\tau<\min\{2+\varepsilon, m(1- \frac{p}{s})\}$
we have
\begin{equation}\label{Morrey-reg-u}
D\,u \,\in L^{p, \tau}(\Omega_0, {\R}^{mn})
\end{equation}
where $\Omega_0$ is a relatively open subset of $\overline{\Omega}$ which
satisfies
$$
\overline{\Omega}\setminus\Omega_0
= \{ x \in \Omega \colon \liminf_{R \to 0} \frac
{1}{R^{m-p}}\int_{\Omega(x,R)} |Du(y)|^p dy > 0
\}.
$$
Moreover, we have
\[
{\cal H}^{m- p -\delta}(\overline{\Omega} \setminus \Omega_0)=0
\]
for some $\delta >0$, where ${\cal H}^r$ denotes the $r$-dimensional
Hausdorff measure.
\end{theorem}

As a corollary of the above theorem we have the following partial
H\"{o}lder regularity result.

\begin{corollary}\label{Cor}
Let
$g$, $u$ and $\Omega_0$ be as in Theorem \ref{mainTh}.
Assume that $p+2\geq m$
and that $s > \max\{m,p\}$.
Then, for some $\alpha \in (0,1)$, we have
\begin{equation}\label{3.2}
u \in C^{0,\alpha}(\Omega_0, {\R}^n).
\end{equation}
\end{corollary}
Moreover, as a corollary of the proof of Theorem \ref{mainTh},
we have the following full-regularity result for the case that
$A$ does not depend on $u$.

\begin{corollary}\label{Cor2}
Assume that $A$ and $g$ satisfy all assumptions of Theorem \ref{mainTh}
and that $A$ does not depend on $u$.
Let $u$ be a minimizer of ${\cal A}$ in the class $X_g$
then
\begin{equation}\label{Morrey-reg-u-2}
D\,u \,\in L^{p, \tau}(\Omega, {\R}^{mn}).
\end{equation}

Moreover, if $p+2\geq m$
and $s> \max\{m,p\}$
, we have full-H\"{o}lder regularity of $u$.
Namely we have
\[
u \in C^{0,\alpha} (\overline{\Omega}, \R^n).
\]
\end{corollary}
\section{Preliminary Lemmas and Proof of the main results}
\setcounter{equation}{0}
\setcounter{theorem}{0}
Throughout the paper we use the following notation:
\begin{eqnarray}
    & & Q^+(x,R)= \{ y\in\R^m~;~
        |y^\alpha - x^\alpha |<R,~\alpha=1,...,m,~
        ~y^m>0  \} ~~\nn \\
    & & ~~~~~~~~~~~~~~~~~~~~~~~\mbox{for}~x\in \R^m\cap \{x~;~x^m=0\},
    ~~R>0,\nn\\
    & & \Omega (x,R) =Q(x,R)\cap \Omega\nn\\
    & & \Gamma(x,R) = Q(x,R)\cap \partial \Omega.\nn
\end{eqnarray}
When the center $x$ is understood,
we sometimes omit the center and write simply $Q(R)$, $Q^+(R)$ etc.
For the sake of simplicity, we always assume that $0<R<1$ in the
following.

We can always reduce locally to the case of flat boundary,
by means of a diffeomorphism which does not change
properties of the functional
assumed in the conditions (A-1)--(A-4).
More precisely, we can choose a positive constant $R_1$ depending only
on $\partial \Omega$ which has the following properties:
\begin{enumerate}
\item A finite number of cubes $\{Q(x,R_1)\}$ centered at $x\in \partial
\Omega$ cover the boundary. Namely,
\[
   \partial \Omega \subset
   \bigcup_{k=1}^{N} Q(x_k, R_1)~~
   x_k \in \partial \Omega~ k=1,...,N,
\]
\item For every $Q(x_k, 2R_1)$, by means of a suitable diffeomorphism,
we can assume that
$x_k=0$ and that
\[
    \begin{array}{l}
    \ds{\Gamma (x_k, 2R_1)
    = Q(0, 2R_1) \cap \partial \Omega \subset \{ x \in \R^m~
    ; x^m=0\},}\\
    \ds{  Q( x_k,2R_1) \cap \Omega
    = Q^+(0, 2R_1) = \{ x\in \R^m~;~
    |x| < 2R_1,~ x^m >0\}.}
    \end{array}
\]
\end{enumerate}

Let us define so-called \it frozen functional. \rm
For some fixed point $x_0 \in \Omega$ and $R>0$
let us define
$A^0 (\xi)$ and $\A^0(u)$
by
\bea
A^0 (\xi)& & \zume \zume \zume = A_R(u_R,\xi) :=
\intmean_{\Omega(x_0,R)} A(y,u_R,\xi) dy, \\
\nn\\
\A^0(u,\Omega(x_0,R))& & \zume \zume \zume := \int_{\Omega(x_0,R)} A^0(Du)dx,
\eea
where
\[
u_R =u_{x_0,R}= \intmean_{\Omega(x_0,R)} u(y) dy.
\]

For weak solutions of the Euler-Lagrange equation of $\A^0$,
we have the following regularity results.

For interior points, we have the following
. (See \cite[Theorem 3.1]{CamAtti}.)

\begin{lemma}\la{ThCamAtti}
Let $u \in H^{1,p}(\Omega, \R^n)\,$ $p \geq 2,$  a solution of the
system
$$
D_\alpha a^\alpha_i (Du)= 0 ~~~~ i=1,...,n
~~~~~\mathrm{in~}~~\Omega,
$$
in the sense that
$$
\int_\Omega a^\alpha_i (Du)D_\alpha \varphi^i dx = 0,\quad
\forall \varphi \in C^\infty_0(\Omega, \R^n)
$$
under the conditions
\begin{enumerate}
\item[(1)]
$a^\alpha_i (0)=0$
\item[(2)]
there exist two constants $\nu>0$ and $M>0$ such that, $\forall x \in \Omega,
\forall
    \xi , \zeta \in \R^{mn}$
    \[
    \| A(\xi) \|
    \leq M \cdot (1+\| \xi \|^2)^{\frac{p - 2}{2}}
    \]
  \[
    A^{\alpha\beta}_{ij} (\xi) \zeta^i_\alpha \zeta^j_\beta \geq
    \nu \cdot  (1+\|\xi\|^2)^{\frac{p - 2}{2}} \|\zeta \|^2
  \]
where $A =( A^{\alpha\beta}_{ij} ) $
and $A^{\alpha\beta}_{ij} (\xi)= \partial a^\alpha_i (\xi) /\partial \xi^j_\beta$.
\end{enumerate}
Then, $\forall Q(\sigma)= Q(x_0, \sigma)\subset\subset \Omega$ and
$\forall t \in (0,1)$
\begin{equation}
\int_{Q(t \sigma)} |Du|^p dx\leq c\cdot t^{\lambda_0} \cdot \int_{Q(\sigma)}
|Du|^p dx, \qquad \lambda_0= \min\{2+\varepsilon_0, m\}
\end{equation}
for some positive constants $\varepsilon_0$ and $c$ which do not
depend on $t,~\sigma$ and
$ x^0.$
\end{lemma}

In the neighborhood of the boundary, by the proof of
\cite[Theorem 7.1]{CamAtti}, we have the following.

\begin{lemma}\la{ThCamAtti-B}
Let $a^\alpha_i (\xi)$ and $\lambda_0$ be as in Lemma \ref{ThCamAtti} and
$v \in H^{1,p}(Q^+(0,R))$ a solution of the problem
\begin{equation}\la{eq-v-B}
    \left\{
\begin{array}{ll}
           \ds{\int_{Q^+(0,R)} a_i^\alpha (Dv+Dg) D_\alpha \varphi^i dx
           =0 } ~ & \forall \varphi \in C^\infty_0 (Q^+(0,R)),\\
           v = 0 & \mbox{on} ~~ \Gamma (0,R),\\
        \end{array}
    \right.
\end{equation}
where $g$ is a given function with
\[
Dg \in L^s(Q^+(0,R))
\]
for some $s>p$.
Then, for every $x_0 \in \Gamma(0,R)$ and
$\tau_0$ with $0<\tau_0< \min\{\lambda_0, m(1-p/s)\}$,
there exist a constant $c>0$ such that
\begin{equation}\la{bdry-est-v}
   \begin{array}{ll}
       & \ds{ \int_{Q^+(x_0,t\sigma)} |W(Dv)|^2 dx }\\
       \leq & \ds{ c  t^{\tau_0}
       \int_{Q^+(x_0, \sigma)}
       |W(Dv)|^2 dx
      +c \sigma^{\tau_0} \big(
       \int_{Q^+(x_0, \sigma)} |W(Dg)|^\frac{2s}{p} dx
       \big)^\frac{p}{s} ,}
   \end{array}
\end{equation}
for any $\sigma \in (0,
R-|x_0|]$ and $t\in (0,1)$,  where
\[
    W(\xi)= ( 1 + |\xi|^2)^{\frac{p-2}{4}}\xi.
\]
\end{lemma}
\noindent 
\it
Outline of the proof.
\rm
Since \pref{eq-v-B} is exactly (7.6) of \cite{CamAtti},
we can proceed as in p.148--150 of \cite{CamAtti}
and get the following estimates:
\[
    \begin{array}{ll}
    & \ds{\int_{Q^+(x_0,t\sigma)} |W(Dv)|^2 dx}\\
    \leq & \ds{c_1
 t^\lambda \int_{Q^+(x_0,\sigma)} |W(Dv)|^2 dx
    + c_1 \int_{Q^+(x_0, \sigma)} (1+|Dv|+|Dg|)^{p-2}|Dg|^2 dx, }\\
    & \\
    & \ds{ \int_{Q^+(x_0, \sigma)} (1+|Dv|+|Dg|)^{p-2}|Dg|^2 dx} \\
    \leq & \ds{ c_2 \int_{Q^+(x_0, \sigma)}|W(Dg)|^2dx
    + c_2 \int_{Q^+(x_0, \sigma)}|Dv|^{p-2}|Dg|^2dx.}\\
       & \\
       & \ds{\int_{Q^+(x_0, \sigma)}|Dv|^{p-2}|Dg|^2dx}\\
       \leq & \ds{\big( 1-\frac{2}{p} \big) \delta
       \int_{Q^+(x_0, \sigma)}
       |W(Dv)|^2 dx + \frac{2}{p} \delta^{1-p/2}
       \int_{Q^+(x_0, \sigma)} |W(Dg)|^2 dx
      }
    \end{array}
 \]
for any $\delta>0$.
These estimates are nothing else than (17)--(19) of
\cite{CamAtti}.
Combining them, we get
\begin{equation}
    \begin{array}{ll}
    & \ds{\int_{Q^+(x_0,t\sigma)} |W(Dv)|^2 dx}\\
    \leq & \ds{c_1 \{  t^\lambda +
    c_2(1-\frac{2}{p})\delta\}
    \int_{Q^+(x_0, \sigma)} |W(Dv)|^2 dx
    }\\
    & \ds{~~~~~~~~~~~~~~~~~~~~~+ c_1 c_2 ( 1+
\frac{2}{p} \delta^{1-p/2}  ) \int_{Q^+(x_0, \sigma)} |W(Dg)|^2 dx.}\\
    \leq & \ds{c_1\{  t^\lambda +c_1 c_2(1-\frac{2}{p})\delta\}
    \int_{Q^+(x_0, \sigma)} |W(Dv)|^2 dx
    }\\
    & \ds{ ~~~~~~~~~~~~~~~~~~~+ c_3 (p,\delta)
    \sigma^{m(1-p/s)}\big( \int_{Q^+(x_0,\sigma)} |W(Dg)|^{2s/p} dx
    \big)^{p/s}.  }
    \end{array}
\end{equation}
Now, using "A useful lemma" of \cite[p.44]{Gia-book2}, we get \pref{bdry-est-v}.
\qed

Moreover, we have the following $L^q$-estimate for $u$.
%
%
\begin{lemma}\label{Lp-Lq-est}
Assume that $u \in H^{1,p}(Q^+(0,R))$ satisfies
\[
   {\cal A}(u, Q^+(0,R)) \leq {\cal A}(u +\varphi,
   Q^+(0,R)) ~~~~\varphi \in
   H^{1,p}_0 (Q^+ (0,R)),
\]
and that $u=g$ on $\Gamma (0,R)$ for some $g \in H^{1,q_1}(Q^+(0,R))$
with $q_1 >p$.
Then there exists an exponent $q \in (p,q_1]$ such that
$u \in H^{1,q}(Q^+(0,r))$ for any $r<R$.
Moreover, if $x_0 \in Q^+(0,r)\cup\Gamma(0,r)$ and $\rho< R-r$,
we have the estimate
\begin{equation}\la{Lp-Lq-est-B}
  \begin{array}{ll}
  & \ds{
  \left(\intmean_{Q(x_0, \rho/2)\cap Q^+(0,R)}
  (1+|Du|^2)^{q/2} dx \right)^{1/q}}\\
  \leq &
  \ds{c \left( \intmean_{Q(x_0, \rho)\cap Q^+(0,R)}
  (1+|Du|^2)^{p/2} dx \right)^{1/p}
  + c   \left(\intmean_{Q(x_0, \rho)\cap Q^+(0,R)} (1+|Dg|^2)^{q/2}
   dx \right)^{1/q}.}
   \end{array}
\end{equation}
In addition, if $Q(x_0,\rho) \subset \subset Q^+(0,R)$, then
we have
\begin{equation}\la{Lp-Lq-est-i} \left(\intmean_{Q(x_0, \rho/2)}
  (1+|Du|^2)^{q/2} dx \right)^{1/q}
  \leq c
  \left( \intmean_{Q(x_0, \rho)}
  (1+|Du|^2)^{p/2} dx \right)^{1/p}.
\end{equation}
\end{lemma}
\noindent 
\it Outline of the Proof.
\rm
For the case that
$Q(x_0,\rho) \subset \subset Q^+(0,R)$, we can proceed as in the proof
of \cite[Theorem 4.1]{G-GActa82} to get \pref{Lp-Lq-est-i}.
For general case, mentioning the difference on the growth conditions,
we can proceed as in the proof of \cite[Lemma 1]{JM}.
\qed

Mention that the above lemma is valid for minimizers of $\A^0$
also.

For bounded domain $D$ with
smooth boundary,
covering $\partial D$ with a finite number of cubes and
using the above local estimates we get the following global
$L^q$-estimates for a minimizer.

\begin{corollary}\la{gl-Lq-est}
Let $D \subset \R^m$ be an open set with
smooth boundary $\partial D$, and let
$v \in H^{1,p}(D)$ be a minimizer for the functional
${\cal A}$ (or $\A^0$) in the class
\[
   X_g := \{ w \in H^{1,p} (D)
   ; w-g \in H^{1,p}_0 (D) \}
\]
for a given map $g \in H^{1,q_1}(D)$, $q_1>p$.
Then $Dv \in L^q (D)$ for some $q \in (p,q_1)$
and
\begin{equation}\la{gl-Lq-est-v}
   \int_D (1+|Dv|^2)^{q/2} dx
  \leq c \int_D
  (1+|Dg|^2)^{q/2} dx.
\end{equation}
\end{corollary}

We show the partial regularity of $u$ by comparing $u$ with $v$.
For this purpose, we need the following lemma which can be shown as
\cite[Theorem 4.2, (4.8) ]{GM-Manus}.
%
\begin{lemma}\la{u-v(lem)}
Let $v \in H^{1,p}(\Omega(x_0,r))$ is a minimizer for
${\cal A}^0 (w, \Omega(x_0, r))$ in the class
\[
     \{ w\in H^{1,p}(\Omega(x_0,r))~;~
     w-u \in  H^{1,p}_0(\Omega(x_0,r))\}
\]
for a given function $u \in H^{1,p}(\Omega(x_0,r))$.
Then we have
\begin{equation}\la{u-v}
     \int_{\Omega(x_0,r)} |Du-Dv|^p dx
     \leq c \big\{ \A^0(u;\Omega(x_0,r)) - \A^0(v;\Omega(x_0,r)) \big\}.
\end{equation}
\end{lemma}

~~

Now, we can prove our main theorem.

\noindent
\it Proof of Theorem \ref{mainTh}. \rm
Assume that $Q(R) = Q(x_0, R) \subset \subset \Omega$.
Let $v \in H^{1,p}(Q(R))$ be a minimizer of $\A^0 (\tilde{v}, Q(R))$
in the class
\[
\{ \tilde{v} \in H^{1,p}(Q(R))~;~ u-\tilde{v}
\in H^{1,p}_0(Q(R))\},
\]
and let $w=u-v$.
First we will estimate $\int_{Q(R)} |Dw|^p dx$.
By Lemma \ref{u-v(lem)} we can see that
\begin{equation}\label{est-Dw}
  \begin{array}{rl}
     & \ds{ \int_{Q(R)} |Dw|^p dx =  c\big\{\A^0(u) - \A^0(v) \big\}}\\
     & \\
     \leq  & \ds{c\int_{Q(R)} \big| A_R(u_R, Du)-A(x,u_R, Du) \big| dx }\\
     & \\
     & \ds{  ~~~~~~~~~~
     + c\int_{Q(R)} \big| A(x, u_R, Du) -A(x,u, Du) \big| dx}\\
     & \\
     & \ds{+ c\int_{Q(R)} \big| A(x,v,Dv) -A(x,u_R, Dv) \big| dx}\\
     & \\
     & \ds{~~~~~~~~~~~~
     +c\int_{Q(R)} \big| A(x,u_R, Dv) -A_R(u_R, Dv) \big| dx}.
   \end{array}
\end{equation}
Here we have used the minimality of $u$.
So, using the assumptions on $A$, we get
\begin{equation}\label{est-Dw(2)}
    \begin{array}{rl}
    & \ds{\int_{Q(R)} |Dw|^p dx}\\
    & \\
    \leq & \ds{ \int_{Q(R)} \big\{ \sigma(x,R) + \omega(|u-u_R|^2) \big\}
    ( 1 + |Du(x)|^2 )^{\frac{p}{2}} dx }\\
    & \\
    & \ds{ + \int_{Q(R)} \big\{ \sigma(x,R) + \omega(|v-u_R|^2) \big\}
    ( 1 + |Dv (x)|^2 )^{\frac{p}{2}} dx. }
    \end{array}
\end{equation}
Using H\"older's inequality,
Lemma \ref{Lp-Lq-est}, \pref{Lp-Lq-est-i}
and the boundedness of $\omega$ and $\sigma$,
we have
\begin{equation}\label{sigmaomegaDu}
\begin{array}{rl}
   & \ds{\int_{Q(R)} \big\{ \sigma(x,R) + \omega(|u-u_R|^2) \big\}
    ( 1 + |Du(x)|^2 )^{\frac{p}{2}} dx}\\
   & \\
   \leq & \ds{ C \left\{ \left(
   \intmean_{Q(R)} \sigma(x,R) dx \right)^{\frac{q-p}{q}}
   + \left(\intmean_{Q(R)} \omega(|u-u_R|^2) dx \right)^{\frac{q-p}{q}}
   \right\} }\\
   & \\
   & \ds{~~~~~~~~~~~~~~~~~~~~~~ \cdot
   \int_{Q(2R)} ( 1 + |Du(x)|^2 )^{\frac{p }{2}} dx.}
\end{array}
\end{equation}
Using Corollary \ref{gl-Lq-est}, and \pref{Lp-Lq-est-i} we get similarly
\begin{equation}\label{sigmaomegaDv}
\begin{array}{rl}
   & \ds{\int_{Q(R)} \big\{ \sigma(x,R) + \omega(|v-u_R|^2) \big\}
     ( 1 + |Dv(x)|^2 )^{\frac{p }{2}} dx}\\
   & \\
   \leq & \ds{ C \left\{ \left(
   \intmean_{Q(R)} \sigma(x,R) dx \right)^{\frac{q-p}{q}}
   + \left(\intmean_{Q(R)} \omega(|v-u_R|^2) dx \right)^{\frac{q-p}{q}}
   \right\} }\\
   & \\
   & \ds{~~~~~~~~~~~~~~~~~~~~~~ \cdot
   \int_{Q(2R)} ( 1 + |Du(x)|^2 )^{\frac{p }{2}} dx.}
\end{array}
\end{equation}
By virtue of concavity of $\omega$, using Jensen's inequality and Poincar\'e
inequality, we have
\begin{equation}\label{est-omega-uv}
    \begin{array}{cl}
      & \ds{ \intmean_{Q(R)} \omega(|u-u_R|^2) dx,~~~~~
      \intmean_{Q(R)} \omega(|v-u_R|^2) dx }\\
      & \\
      \leq \!\!\! & \ds{
      C\omega \left( R^{p-m}\int_{Q(R)} |Du|^p dx \right).}
    \end{array}
\end{equation}
Combining \pref{est-Dw(2)} -- \pref{est-omega-uv}, we obtain
\begin{eqnarray}
     &  & \int_{Q(R)} |Dw|^p dx \nonumber\\
     & & \nonumber \\
     & \leq & \!\!\! C \left\{
     \left(\intmean_{Q(R)} \sigma(x,R) dx \right)^{\frac{q-p}{q}}
     + \omega \left( R^{p - m}\int_{Q(R)} |Du|^p dx \right)^{\frac{q-p}{q}}
     \right\} \\
     & & \nonumber \\
     & & ~~~~~~~~~~~~~~~~~~~~~~~~~~~\cdot \int_{Q(2R)} ( 1 + |Du(x)|^2 )^{\frac{p }{2}} dx.\nonumber
     \label{est-Dw(3)}
\end{eqnarray}
Now, from Lemma \ref{ThCamAtti} and
the above inequality, we get
\begin{eqnarray}
    & & \int_{Q(r)} |Du|^p dx
        \leq \int_{Q(r)}\big(|Dv|^p + |Dw|^p\big )dx \nonumber\\
    & & \nonumber\\
    & \quad \qquad \leq & \!\!\! C \left\{ \left( \frac{r}{R} \right)^\lambda
    +\left(\intmean_{Q(R)} \sigma(x,R) dx \right)^{\frac{q-p}{q}}
    \right. \label{est-Du-final}\\
    & & \nonumber\\
    & & ~~~\left.
        + \omega \left( R^{p - m }\int_{Q(R)} |Du|^2 dx \right)^{\frac{q-p}{q}}
        \right\} \cdot \int_{Q(2R)} ( 1 +
Du(x)|^p)^{\frac{p}{2}}
 dx. \qquad\qquad
    \nonumber
\end{eqnarray}

Let us consider the behavior of $u$ near the boundary.
Let $Q(x_l, 2R_1)$ be a member of the covering
$\{ Q(x_k,2R_1) \}$ which is introduced at the beginning of
this section.
Then, $u$ satisfies
\begin{equation}\la{u-bdry}
\left\{
    \begin{array}{l}
    \ds{\A(u, Q^+(x_l,2R_1)) \leq \A(u+\varphi,Q^+(x_l,2R_1)) ~~
    \forall \varphi \in H^{1,p}_0(Q^+(x_l,2R_1)),}\\
    u=g ~~ \mbox{on}~~\Gamma(x_l,2R_1) .
    \end{array}
\right.
\end{equation}
Fix a point $x_0\in \Gamma(x_l,R_1)$ and a positive number
$R< R_1$ arbitrarily
(here, mention that $Q^+(x_0,R)
\subset Q^+(x_l,2R_1)$).
Let $v \in H^{1,p} (Q^+(x_0,R))$ be a minimizer of
$\A^0(v,Q^+(x_0,R))$ in the class
\[
\{ v \in H^{1,p}(Q^+(x_0,R))~;~ u-v\in H^{1,p}_0(Q^+(x_0,R))\},
\]
and put $w=u-v$.
Then, using Lemma \ref{u-v(lem)}, we can proceed as in the interior case
and get
\begin{equation}\label{est-Dw(2)-B}
    \begin{array}{rl}
    & \ds{\int_{Q^+(R)} |Dw|^p dx}\\
    & \\
    \leq & \ds{ \int_{Q^+(R)} \big\{ \sigma(x,R) + \omega(|u-u_R|^2) \big\}
    ( 1 + |Du(x)|^2 )^{\frac{p}{2}} dx }\\
    & \\
    & \ds{ + \int_{Q^+(R)} \big\{ \sigma(x,R) + \omega(|v-u_R|^2) \big\}
    ( 1 + |Dv(x)|^2 )^{\frac{p}{2}} dx. }
    \end{array}
\end{equation}
Moreover, using \pref{Lp-Lq-est-B} instead of
\pref{Lp-Lq-est-i}
and proceeding
as in the interior case, we have
\begin{eqnarray}
     &  & \int_{Q^+(R)} |Dw|^p dx \nonumber\\
     & & \nonumber \\
     & \leq & \!\!\! C \left\{
     \left(\intmean_{Q^+(R)} \sigma(x,R) dx \right)^{\frac{q-p}{q}}
     + \omega \left( R^{p - m}\int_{Q^+(R)} |Du|^p dx \right)^{\frac{q-p}{q}}
     \right\} \nn \\
     & & \label{est-Dw(3)-B}\\
     & & ~~~~~~~~~~~~~~~~~~~~~~~~~~~
     \cdot \int_{Q^+(2R)} ( 1 + |Du(x)|^2 )^{\frac{p }{2}} dx
     \nonumber\\
     & & \nn\\
     & & ~~~~ + CR^{m\frac{q-p}{q}}\left(
     \int_{Q^+(2R)} (1+|Dg|^2)^\frac{q}{2} dx\right)^\frac{p}{q}.\nn
\end{eqnarray}
Now, combining \pref{bdry-est-v}
and \pref{est-Dw(3)-B}, we obtain
\begin{eqnarray}
    & & \int_{Q^+(r)} |Du|^p dx \nonumber\\
    & & \nonumber\\
    & \quad \qquad \leq & \!\!\! C \left\{
    \left( \frac{r}{R} \right)^{\tau_0}
    +\left(\intmean_{Q^+(R)} \sigma(x,R) dx \right)^{\frac{q-p}{q}}
    + \omega \left( R^{p - m }\int_{Q^+(R)}
    |Du|^p dx \right)^{\frac{q-p}{q}}
        \right\}
    \nn\\
    & & \label{est-Du-semifinal-B}\\
    & & ~~~~~~~~~~~~~~\cdot \int_{Q^+(2R)} ( 1 + |Du(x)|^2 )^{\frac{p}{2}}
    dx \qquad\qquad \nn\\
    & & + c  R^{\tau_0} \left(
    \int_{Q^+(R)} (1+|Dg|^2)^\frac{s}{2} dx \right)^\frac{p}{s}
    + CR^{m\frac{q-p}{q}}\left(
     \int_{Q^+(2R)} (1+|Dg|^2)^\frac{q}{2} dx\right)^\frac{p}{q}.\nn
    \nonumber
\end{eqnarray}
Since we are assuming that $Dg \in L^s$ for some $s>p$,
and we can choose $q>p$ sufficiently near to $p$, without
loss of generality we can assume that $s>q>p$.
So, we can estimate the last term of \pref{est-Du-semifinal-B}
as follows:
\begin{eqnarray}
R^{m\frac{q-p}{q}}\left(
     \int_{Q^+(2R)} (1+|Dg|^2)^\frac{q}{2} dx\right)^\frac{p}{q}
& \leq & CR^{m(1-p/s)}
\left(\int_{Q^+(2R)}(|1+|Dg|^2)^{\frac{s}{2}}dx \right)^\frac{p}{s}.   \nn
\end{eqnarray}
Here, we can assume that $R<1$, so the above estimates hold
even if $m(1-p/s)$ can be replaced by the smaller constant $\tau_0$.
Mentioning the above fact and
combining the above estimate with \pref{est-Du-semifinal-B}, we
get the following estimate.
\begin{eqnarray}
    & & \int_{Q^+(r)} |Du|^p dx \nonumber\\
    & & \nonumber\\
    &
    \leq & \!\!\! C
    \left\{ \left( \frac{r}{R} \right)^{\tau_0}
    +\left(\intmean_{Q^+(R)} \sigma(x,R) dx \right)^{\frac{q-p}{q}}
    + \omega \left( R^{p - m }\int_{Q^+(R)}
    |Du|^p dx \right)^{\frac{q-p}{q}}
        \right\}
    \label{est-Du-final-B}\\
    & & \nn\\
    & & ~~~~~~~~~~~~~~\cdot \int_{Q^+(2R)} ( 1 + |Du(x)|^2 )^{\frac{p}{2}}
    dx + C(g)R^{\tau_0}. \nn
\end{eqnarray}

By the assumption (A-1), we have
\[
\intmean_{Q(R)} \sigma(x,R) dx \to 0~~~~\mathrm{as}~~R\to 0.
\]
So, using "A useful Lemma" on p.44 of \cite{Gia-book2} for
\pref{est-Du-final} and \pref{est-Du-final-B}, and
putting
\begin{equation}
   \Phi(x,r) = \int_{\Omega(x,r)} (1+|Du|^2)^\frac{p}{2} dx,
\end{equation}
we can see that
for any $\tau$ with $ 0 < \tau < \tau_0(<\lambda_0)$
there exist positive constants $\delta$, $M$ and
$R_0~(R_0<R_1/2)$ with the following properties.
\begin{description}
\item[{[}Interior Case{]}]
If
\[
r_1,~~~ r_1^{p-m} \Phi (x,r_1) < \delta
\]
for some $r_1 \in (0, R_0) $ with $Q(x,r_1) \subset \subset \Omega$,
then for $0 < \rho < r < r_1$ we have
\begin{equation}\la{Morrey-est1-i}
    \Phi(x,\rho) \leq M\left( \frac{\rho}{r}\right)^\tau
    \Phi(x,r).
\end{equation}
\item[{[}Boundary Case{]}]
For $x \in \partial \Omega$, if
\[
r_1,~~~r_1^{p-m} \Phi (x,r_1) < \delta
\]
for some $r_1\in (0, R_0) $, then we have
\begin{equation}\la{Morrey-est1-b}
    \Phi(x,\rho) \leq M\left( \frac{\rho}{r}\right)^\tau
    \Phi(x,r)+ M \rho^\tau.
\end{equation}
\end{description}

Now, we can proceed as in pp.318--319 of Giusti's book
\cite{Giusti-book} to show partial Morrey-type regularity of $u$.
Namely, there exist positive constants $\delta$ and $M$ with the
following properties. For any $x \in \Omega$, if
\begin{equation}
   r_0, ~~r_0^{p-m} \Phi (x, r_0) \leq \delta
\end{equation}
for some $r_0>0$, then
\begin{equation}
\rho^{-\tau} \Phi (x,\rho) \leq \tilde{M}.
\end{equation}
So, we get the assertion.
\qed

~~

\noindent
\it Proof of Corollary \ref{Cor}. \rm
When $p+2\geq m$ and $s>\max\{ m,p\}$, we can
take $\tau$ sufficiently near to
$\min\{2+\varepsilon, m(1-\frac{p}{s})\}$
so that $\tau > m-p$.  So,
Corollary \ref{Cor} is a direct consequence of Theorem \ref{mainTh}
and Morrey's theorem on the growth of the Dirichlet integral
(see, for example, p.43 of \cite{Gia-book2}).
\qed

~~

\noindent
\it
Proof of Corollary \ref{Cor2}.
\rm When $A(x,u,\xi)$ does not depend on $u$, we can proceed
as in the proof of Theorem \ref{mainTh} without the term with
$\omega$ and get, instead of
\pref{est-Du-final} and \pref{est-Du-semifinal-B},
\begin{eqnarray}
    & & \int_{Q(x_0, r)} |Du|^p dx
        \nonumber\\
    & & \label{est-Du-final-2}\\
    & \quad \qquad \leq & \!\!\! C \left\{ \left( \frac{r}{R} \right)^\lambda
    +\left(\intmean_{Q(R)} \sigma(x,R) dx \right)^{\frac{q-p}{q}}
        \right\} \cdot \int_{Q(2R)} ( 1 + |Du(x)|^2 )^{\frac{p}{2}}
 dx. \qquad\qquad
    \nonumber
\end{eqnarray}
for $ Q(2R)=Q(x_0,2R)\subset\subset \Omega$ and
\begin{eqnarray}
    & & \int_{Q^+(x_0,r)} |Du|^p dx \nonumber\\
    & & \label{est-Du-final-B-2}\\
    & \quad \qquad \leq & \!\!\! C \left\{ \left( \frac{r}{R} \right)^\lambda
    +\left(\intmean_{Q^+(R)} \sigma(x,R) dx \right)^{\frac{q-p}{q}}
        \right\}\int_{Q^+(2R)} ( 1 + |Du(x)|^2 )^{\frac{p}{2}}
    dx + C(g)R^\tau,\nn
\end{eqnarray}
for $x_0\in \partial \Omega$.
So, we can proceed as in the last part of Theorem \ref{mainTh}
without assuming that
\[
r_1^{p-m} \Phi(x,r_1) = r_1^{p-m}
\int_{\Omega(x,r_1)} (1+|Du|^2)^\frac{p}{2} dx < \delta.
\]
and see that
\[
\rho^{-\tau} \Phi (x,\rho) \leq \tilde{M}
\]
for all $x\in \Omega$. Thus we get the assertions
\qed

\begin{remark}\label{4-2}
Without any restriction on the dimension of the domain,
it is not possible to obtain H\"older regularity result in all the
domain $\Omega$ as showed by V.
\v{S}verak and X. Yan in a counterexample contained in \cite{SY}.
\end{remark}

\end{document}